\documentclass{amsart}[12pt]
\usepackage{mathrsfs, amsmath,amssymb}
\usepackage{fullpage}
\newtheorem{teo}{Theorem}
\newtheorem{pro}{Proposition}
\newtheorem{lem}{Lemma}
\newtheorem{cor}{Corollary}
\newtheorem{rem}{Remark}
\newtheorem*{rems}{Remarks}

\author{L. Deleaval}
\address{Laboratoire d'Analyse et de Math\'ematiques appliqu\'ees \\ Universit\'e Paris-Est Marne-la-Vall\'ee \\
France}
\email{luc.deleaval@u-pem.fr}

\author{N. Demni}
\address{IRMAR, Universit\'e de Rennes 1\\ Campus de
Beaulieu\\ 35042 Rennes cedex\\ France}
\email{nizar.demni@univ-rennes1.fr}
\subjclass[2010]{33C45; 33C52; 33C65}
\keywords{Modified Bessel functions; Gegenbauer polynomials; Generalized Bessel function; Dihedral groups; Elementary symmetric functions.}
\title{On a Neumann-type series for modified Bessel functions of the first kind}
\begin{document}
\maketitle

\begin{abstract}
In this paper, we are interested in a Neumann-type series for modified Bessel functions of the first kind which arises in the study of Dunkl operators associated with dihedral groups and as an instance of the Laguerre semigroup constructed by Ben Said-Kobayashi-Orsted. We first revisit the particular case corresponding to the group of square-preserving symmetries for which we give two new and different proofs other than the existing ones. The first proof uses the expansion of powers in a Neumann series of Bessel functions while the second one is based on a quadratic transformation for the Gauss hypergeometric function and opens the way to derive further expressions when the orders of the underlying dihedral groups are powers of two. More generally, we give another proof of De Bie \& al formula expressing this series as a $\Phi_2$-Horn confluent hypergeometric function. In the course of proving, we shed the light on the occurrence of multiple angles in their formula through elementary symmetric functions, and get a new representation of Gegenbauer polynomials.  
\end{abstract}

\section{Introduction}
In this paper, we are interested in the following series:
\begin{equation}\label{Series}
F_{n,k}(R,\xi) := \left(\frac{2}{R}\right)^{nk} \sum_{j \geq 0}(j+k){\it I}_{n(j+k)}(R) C_j^{(k)}(\cos \xi).
\end{equation}
Here, $R, k$ are positive reals, $\xi \in [0,\pi]$, $n \geq 1$ is an integer, and ${\it I}_{\nu}, C_j^{(k)},$ denote respectively the modified Bessel function of the first kind of index $\nu > -1/2$ and the $j$-th Gegenbauer polynomial of parameter $k$ (see for instance \cite[page 222]{AAR} and \cite[page 302]{AAR} for their respective definitions). Moreover, $F_{n,k}$ is well-defined since the following estimate holds (\cite{AAR}, p. 302),: 
\begin{equation}\label{Estimate}
|C_j^{(k)}(\cos \xi)| \leq C_j^{(k)}(1) = \frac{(2k)_j}{j!}.
\end{equation}
We shall call the series defining $F_{n,k}(R,\xi)$ `of Neumann-type' since it is a Neumann series for modified Bessel functions when $n=1$ (\cite[chapter XVI]{Wat}). When $n=2$, it appears in \cite{Kob-Man} in relation with certain representations emanating from conformal geometry considerations of the indefinite orthogonal group of rank two  and in \cite{BKO} as the kernel of the so-called Laguerre semi-group where the parameter $n$ is allowed to be positive real. In an apparently different setting, it encodes the so-called generalized Bessel function associated with dihedral groups, see \cite{Demni1}, where it was given a variational formula provided $k \geq 1$ is an integer. In particular, when $n=1$, the series \eqref{Series} reduces to an exponential function (see for instance \cite[page 369]{Wat}): 
\begin{equation}\label{GegenExp}
F_{1,k}(R,\xi) = \frac{1}{\Gamma(k)}e^{R\cos \xi},
\end{equation}
a known result going back to L. Gegenbauer and often helpful to compute Euclidean Fourier transforms of zonal functions. Besides, if $n=2$, then it was proved in \cite{Kob-Man} and independently in \cite{Demni1} and \cite{DeB} that $F_{2,k}(R,\xi)$ may be expressed through the normalized modified Bessel function of the first kind in the variable $\cos(\xi/2)$: 
\begin{teo}[\cite{Demni1,Kob-Man}] \label{Th1}
Let $k, R$ be two positive real numbers and let $\xi \in [0,\pi]$. Then,
\begin{equation}\label{Bessel}
F_{2,k}(R, \xi) = \frac{1}{2\Gamma(2k)} \mathcal{I}_{k-(1/2)}\left(R\cos(\xi/2)\right).
\end{equation}
\end{teo}
In \eqref{Bessel}, we denoted by
\begin{equation*}
\mathcal{I}_{k-(1/2)}(u) := \sum_{j \geq 0} \frac{1}{j!(k+(1/2))_j} \left(\frac{u}{2}\right)^{2j} = \Gamma\left(k+\frac{1}{2}\right) \left(\frac{2}{u}\right)^{k-(1/2)} {\it I}_{k-(1/2)}(u)
\end{equation*}
the normalized modified Bessel function of the first kind. 

More generally, the variational formula proved in \cite{Demni1} suggests that for any $n \geq 1$ and any positive real $k$, $F_{n,k}(R,\cdot)$ is a function of the angles 
\begin{equation}\label{angles}
\theta_s := \frac{\xi+2\pi s}{n}, \quad s = 1,\ldots, n, 
\end{equation}
which correspond from a geometrical point of view to the $n$-fold covering of the unit circle $z \mapsto z^n$. Recently, this suggestion was confirmed in \cite{CDBL} where \eqref{Series} is expressed by means of Horn's hypergeometric function $\Phi_2$ (see the definition below) in the variables 
\begin{equation*}
\cos(\theta_s), \quad s=1,\ldots, n.
\end{equation*}
The approach undertaken in \cite{CDBL} is based on introducing a new parameter in \eqref{Series} with respect to which the authors compute the Laplace transform of the new series and put the obtained expression in a suitable (product) form allowing its inversion. However, the occurrence of the angles $\theta_s$ looks rather mysterious in the sense that it is not obvious how the Horn series $\Phi_2$ in the variables $\cos(\theta_s), s = 1,\ldots, n,$ sum up to \eqref{Series} which only depends on $\cos \xi$. Looking for an explanation of this occurrence is the starting point of our investigations below where we shall justify it by proving that all the elementary symmetric functions in these variables do not depend on $\xi$ except the $n$-th one. Our proof uses induction as well as Newton's identities relating elementary symmetric functions and power sums (\cite{MD}). We shall also revisit the De Bie \& al. formula which turns out to be equivalent to an interesting identity extending the known `half-angle' expansion of $\cos(\xi/2)$ as a finite sum of Gegenbauer polynomials in the variable $\cos \xi$ (\cite{Erd1}). Actually, we write another proof of this identity where the polynomial whose coefficients are the elementary symmetric functions alluded to above is identified through a reverse Tchebycheff polynomial. Furthermore, we obtain from the same identity, by inverting a triangular system. a new representation of Gegenbauer polynomials as a finite sum in the variables $\cos(\theta_s), s = 1,\ldots, n,$. In this respect, we point out that a similar, yet different, representation already appeared in the study of generalized Fibonacci polynomials and is only valid when the parameter $k$ is a positive integer (\cite{Dil}). 

Before dealing with the series \eqref{Series} for arbitrary $n$, we give again some interest in the series \eqref{Series} with $n=2$ corresponding to the dihedral group preserving the square and supply two new and different proofs other than the existing ones (\cite{Kob-Man}, \cite{Demni1}, \cite{DeB}). The first proof uses a series arrangement technique together with the expansion of powers as Neumann series of modified Bessel functions of the first kind (see for instance \cite[page 138]{Wat}). As to the second one, it appeals to an interesting integral representation of Gegenbauer polynomials which in some sense doubles both their degrees and parameters and rather bisects the angles of their arguments. While this integral representation was noticed in \cite{Xu} as a consequence of Dijksma-Koornwinder product formula for Jacobi polynomials (\cite{Dij-Koo}), we notice here that it also follows from a quadratic transformation of the Gauss hypergeometric function. The reason behind doing so is to wonder whether or not analogous integral representations holds for larger values of $n \geq 3$. Moreover, it is worth noting that successive applications of the aforementioned integral representation provides further expressions of \eqref{Series} corresponding to dihedral groups whose orders are powers of two. 

The paper is organized as follows. The next section contain both proofs of Theorem \ref{Th1} below. In section 3, we recall De Bie \& al. formula and deduce from it a miscellaneous identity satisfied by Gegenbauer polynomials for which we give another proof. 
In the same section, we prove by induction that all the elementary symmetric functions of orders less than $n-1$ in the variables $\cos(\theta_s), s=1,\ldots, n,$ do not depend on $\cos \xi$ and give their explicit expressions using reversed Tchebycheff polynomials.   
The paper is closed with a new representation of Gegenbauer polynomials through the variables $\cos(\theta_s), s=1,\ldots, n$. 

\section{The case $n=2$: yet two other proofs}
In this section, we give two new and different proofs of \eqref{Bessel}. 

\subsection{\bf First proof: series arrangement}
The first proof of \eqref{Bessel} we supply consists in using the very definition of Gegenbauer polynomials (\cite[table 8.93]{Gra-Ryz}): for every non-negative integer $j$ and every $x \in [-1,1]$, 
\begin{align}\label{ExpGeg}
C_j^{(k)}(x) = (-1)^jC_j^{(k)}(-x) & = (-1)^j \frac{(2k)_j}{j!} \,{}_2F_1\left(-j, j+2k, k+\frac{1}{2}; \frac{1+x}{2}\right) \nonumber
\\& = \sum_{m=0}^j \frac{(-1)^{m+j}}{(j-m)!m!} \frac{(2k+j)_m(2k)_j}{(k+(1/2))_m} \left(\frac{1+x}{2}\right)^m,
\end{align} 
where ${}_2F_1$ stands for the Gauss hypergeometric function (\cite[chapter 2]{AAR}). Summing over $j \geq 0$ and inverting the summation order, we get
\begin{multline*}
F_{2,k}(R, \xi) =  \left(\frac{2}{R}\right)^{2k} \sum_{m \geq 0} \frac{1}{m!(k+(1/2))m}\bigl(\cos(\xi/2)\bigr)^{2m}\sum_{j \geq m}\frac{(-1)^{m+j}}{(j-m)!}(2k+j)_m(2k)_j (j+k){\it I}_{2(j+k)}(R).
\end{multline*}
Using the index change $j \mapsto j-m$ in the inner series and invoking the identity
\[
(2k+j+m)_m(2k)_{j+m}=\frac{\Gamma(2k+j+2m)}{\Gamma(2k+j+m)}\frac{\Gamma(2k+j+m)}{\Gamma(2k)}=(2k)_{j+2m}
\]
lead to
\begin{multline*}
F_{2,k}(R, \xi) = \left(\frac{2}{R}\right)^{2k} \sum_{m \geq 0} \frac{1}{m!(k+(1/2))m}\bigl(\cos(\xi/2)\bigr)^{2m}\sum_{j \geq 0}\frac{(-1)^{j}}{j!} (2k)_{j+2m}(j+m+k){\it I}_{2(j+m+k)}(R). 
\end{multline*}
But the following Neumann series 
\begin{equation*}
\left(\frac{R}{2}\right)^{\nu} = \sum_{j \geq 0} \frac{(-1)^{j}(\nu+2j)\Gamma(\nu+j)}{j!}{\it I}_{\nu+2j}(R), \quad \nu \geq 0,
\end{equation*}
which may be derived from formula (1) in \cite[page 138]{Wat} using a standard variable change, shows that 
\begin{align*}
\sum_{j \geq 0}\frac{(-1)^{j}}{j!} (2k)_{j+2m}(j+m+k){\it I}_{2(j+m+k)}(R) = \frac{1}{2\Gamma(2k)}\left(\frac{R}{2}\right)^{2m+2k}, 
\end{align*}
whence we finally obtain
\begin{align*}
F_{2,k}(R, \xi) = \frac{1}{2\Gamma(2k)}\sum_{m \geq 0} \frac{1}{m!(k+(1/2))m}\bigl(\cos(\xi/2)\bigr)^{2m}\left(\frac{R}{2}\right)^{2m}  = \frac{1}{2\Gamma(2k)}\mathcal{I}_{k-(1/2)}\bigl(R\cos(\xi/2)\bigr).
\end{align*}

\subsection{\bf Second proof: a quadratic transformation for the Gauss hypergeometric function}

The second proof of \eqref{Bessel} we write here is mainly based on the following integral representation, noticed in \cite{Xu} as a particular occurrence of the product formula for Jacobi polynomials due to Dijksma-Koornwinder (\cite{Dij-Koo}): 
\begin{equation}\label{IR1}
C_j^{(k)}(2x^2 -1) = \frac{\Gamma(k+(1/2))}{\sqrt{\pi}\,\Gamma(k)}\int_{-1}^1 C_{2j}^{(2k)}(xz)(1-z^2)^{k-1}dz.
\end{equation}
Indeed, inserting \eqref{IR1} in the series defining $F_{2,k}(R, \xi)$, we get: 
\begin{equation*}
F_{2,k}(R, \xi) = \left(\frac{2}{R}\right)^{2k} \frac{\Gamma(k+(1/2))}{2\sqrt{\pi}\,\Gamma(k)}\int_{-1}^1 \sum_{j \geq 0}(2j+2k) I_{2j+2k}(R) C_{2j}^{(2k)}\bigl(\cos(\xi/2)z\bigr)(1-z^2)^{k-1}dz,
\end{equation*}
where the interchange of the infinite summation and the integral orders is ensured by \eqref{Estimate}. Besides, the symmetry relation $C_j^{(k)}(x) = (-1)^jC_j^{(k)}(-x)$ (\cite[table 8.93]{Gra-Ryz}) entails
\begin{multline*}
F_{2,k}(R, \xi) =  \left(\frac{2}{R}\right)^{2k} \frac{\Gamma(k+(1/2))}{4\sqrt{\pi}\,\Gamma(k)}\int_{-1}^1 \sum_{j \geq 0}(j+2k) I_{j+2k}(R)   \Bigl(C_{j}^{(2k)}\bigl(\cos(\xi/2)z\bigr)+C_{j}^{(2k)}\bigl(-\cos(\xi/2)z\bigr)\Bigr) \\ (1-z^2)^{k-1}dz
   =  \frac{\Gamma(k+(1/2))}{2\sqrt{\pi}\,\Gamma(k)} \\ \int_{-1}^1 \left\{\left(\frac{2}{R}\right)^{2k} \sum_{j \geq 0}(j+2k) I_{j+2k}(R) C_{j}^{(2k)}\bigl(\cos(\xi/2)z\bigr)\right\} (1-z^2)^{k-1}dz.
\end{multline*}
The series between brackets is, up to the subtitution $k \rightarrow 2k$, the Gegenbauer expansion \eqref{GegenExp}. Consequently,
\begin{align*}
F_{2,k}(R, \xi) = \frac{\Gamma(k+(1/2))}{2\sqrt{\pi}\,\Gamma(k)\,\Gamma(2k)}\int_{-1}^1 e^{R\cos(\xi/2)z}(1-z^2)^{k-1}dz
 = \frac{1}{2\Gamma(2k)} \mathcal{I}_{k-(1/2)}\left(R\cos(\xi/2)\right)
\end{align*}
where the last equality follows from the Poisson integral representation of the modified Bessel function $\mathcal{I}_{k-1/2}$ (\cite[page 79]{Wat}). 

\begin{rem} 
The identity \eqref{IR1} may be also derived from a quadratic transformation for the Gauss hypergeometric function which gives the following representation for even Gegenbauer polynomials (see \cite[page 176]{Erd2}):
\begin{align*}
C_{2j}^{(2\nu)}(x) &= (-1)^j\frac{(2\nu)_j}{j!}\,{}_2F_1\left(-j, j+2\nu, \frac{1}{2}; x^2\right).
\end{align*}
Accordingly, we have
\begin{multline*}
\frac{\Gamma(k+(1/2))}{\sqrt{\pi}\,\Gamma(k)}\int_{-1}^1 C_{2j}^{(2k)}(xz)(1-z^2)^{k-1}dz= \\ (-1)^j\frac{(2k)_j}{j!}\frac{\Gamma(k+(1/2))}{\sqrt{\pi}\,\Gamma(k)}\sum_{m=0}^j\frac{(-j)_m}{m!} \frac{(j+2k)_m}{(1/2)_m}x^{2m} \int_{-1}^1 z^{2m}(1-z^2)^{k-1}dz.
\end{multline*}
But, since
\[
\frac{\Gamma(k+(1/2))}{\sqrt{\pi}\,\Gamma(k)}\int_{-1}^1 z^{2m}(1-z^2)^{k-1}dz = \frac{(1/2)_m}{(k+(1/2))_m},
\]
then we obtain 
\begin{align*}
\frac{\Gamma(k+(1/2))}{\sqrt{\pi}\,\Gamma(k)}\int_{-1}^1 C_{2j}^{(2k)}(xz)(1-z^2)^{k-1}dz&=(-1)^j\frac{(2k)_j}{j!}\sum_{m=0}^j\frac{(-j)_m}{m!} \frac{(j+2k)_m}{(k+(1/2))_m}x^{2m}\\
&=C_j^{(k)}(2x^2-1)
\end{align*}
as required. It is then natural to wonder whether an analog of \eqref{IR1} holds for arbitrary $n \geq 2$: are there (non necessarily orthogonal) polynomials $(Q_j^{n,k})_{j \geq 0}$ and a finite measure $\mu^{(k,n)}$ such that 
\begin{equation*}
C_j^{(k)}(\cos \xi) = \int_{\mathbb{R}} Q_{nj}^{(n,k)}(z\cos(\xi/n)) \mu^{(k,n)}(dz)?
\end{equation*} 
\end{rem}

\subsection{Further expressions}
In the previous paragraph, the expression of the series corresponding to $n=2$ is derived from Gegenbauer's expansion ($n=1$) by means of \eqref{IR1}. We can iterate this procedure in order to derive further expressions of \eqref{Series} for integers 
$n = 2^m, m \geq 1$ which correspond to dihedral groups of orders $4n = 2^{m+2}$. For instance, for $n=4$, the same lines of the second proof above yield: 
\begin{multline*} 
\left(\frac{2}{R}\right)^{4k} \sum_{j \geq 0}(j+k) I_{4(j+k)}(R)C_{j}^{(k)}(\cos \xi) = \\ \frac{\Gamma(k+(1/2))}{2\sqrt{\pi}\,\Gamma(k)} \left(\frac{2}{R}\right)^{4k} 
 \int_{-1}^1 \sum_{j \geq 0}(j+2k) I_{2(j+2k)}(R)C_{j}^{(2k)}(z \cos \xi/2)(1-z^2)^{k-1}dz,
\end{multline*}
and then
\begin{multline*} 
\left(\frac{2}{R}\right)^{4k} \sum_{j \geq 0}(j+k) I_{4(j+k)}(R)C_{j}^{(k)}(\cos \xi) = \\ \frac{\Gamma(k+(1/2))}{4\sqrt{\pi}\,\Gamma(k) \Gamma(4k)} 
 \int_{-1}^1 \mathcal{I}_{2k-1/2}\left(R\sqrt{\frac{1+z\cos(\xi/2)}{2}}\right) (1-z^2)^{k-1}dz.
\end{multline*}
The last integral may be written as a hypergeometric series in two variables and we shall not do this here since a more general formula we recall below and valid for arbitrary integers $n \geq 1$ shows that \eqref{Series} may be expressed through one of the numerous Horn's hypergeometric functions (\cite{Erd1}). In the next section, we give another proof of this formula and explain the occurrence there of the multiple angles $\theta_s, s =1, \ldots, n,$ by proving that the elementary symmetric functions $(e_m)_{m=1}^n$ in the variables $\cos(\theta_s), s =1, \ldots, n,$ do not depend on $\xi$ unless $m=n$.

\section{Another proof of De Bie \& al. formula} 
In this section, we assume $n \geq 1$ is an arbitrary but fixed integer and write another proof of De Bie \& al. formula we recall below. We consider, for fixed $R > 0$ and $\xi \in [0,\pi]$, the following function\begin{equation*}
t \in (0,+\infty)\mapsto \left(\frac{2}{R}\right)^{nk} \sum_{j \geq 0}(j+k){\it I}_{n(j+k)}(R t) C_j^{(k)}(\cos \xi).
\end{equation*}
If we set $a_s := R\cos \theta_s$, where $\theta_s$ is defined in \eqref{angles}, then it is proved in \cite{CDBL} that its Laplace transform in the variable $z$ is given by
\begin{equation}\label{LapTra}
-\Gamma(nk) \frac{d}{ds} \frac{1}{(z-a_1)^k\ldots (z-a_n)^k},
\end{equation}
provided $\Re(z) > \max(a_1, \ldots, a_n)$. Viewing this Laplace transform as a function in the shifted variable $z-a_1$ and using formula  (5) in \cite[table 4.24]{Erd},  it was argued in \cite{CDBL} that\footnote{The authors only displayed the formula for $n=3$ but it still holds for arbitrary integers $n \geq 1$ if we assume that the confluent function reduces to one when $n=1$.}:
\begin{equation*}
F_{n.k}(R, \xi) \propto e^{a_1}\Phi^{(n-1)}_{2}(\underbrace{k, \ldots, k}_{n-1}; nk; a_2-a_1, \ldots, a_n-a_1),
\end{equation*}
where for $m \geq 1$, 
\begin{equation*}
\Phi_{2}^{(m)}(\beta_1, \ldots, \beta_{m}; \gamma; x_1, \dots, x_{m}) := \sum_{j_1, \ldots, j_{m} \geq 0} \frac{(\beta_1)_{j_1}\ldots (\beta_{m})_{j_{m}}}{(\gamma)_{j_1+\cdots+j_{m}}} \frac{x_1^{j_1}}{j_1!}\cdots  \frac{x_{m}^{j_{m}}}{j_{m}!}
\end{equation*}
is the Horn confluent hypergeometric function of $m$ variables (see \cite[chapter V]{Erd1}  for more details). However, the same formula (5) in \cite[table 4.24]{Erd} yields 
\begin{align*}
F_{n,k}(R, \xi) \propto  \Phi_{2}^{(n)}(\underbrace{k, \dots, k}_{n}; nk; a_1, a_2, \dots, a_n),
\end{align*}
which on the one hand gives an argument transformation for $\Phi_2^{(m)}$ and on the other hand, is well suited for our purposes. Indeed, the last equality may be further expanded as 
\begin{align*}
\sum_{m,j \geq 0}\frac{n(j+k)\Gamma(nk)}{m!\Gamma(n(j+k)+m+1)}C_j^{(k)}(\cos \xi) \left(\frac{R}{2}\right)^{nj+2m} = \sum_{j_1, \ldots, j_{n} \geq 0} \frac{(k)_{j_1}\ldots (k)_{j_{n}}}{(nk)_{j_1+\cdots+j_{n}}} \frac{a_1^{j_1}}{j_1!}\cdots  \frac{a_{n}^{j_{n}}}{j_{n}!},
\end{align*}
and as such, we end up with the following interesting and miscellaneous identity, in which we set: 
\begin{equation*}
b_s := \cos\theta_s=\cos\left(\frac{\xi +2\pi s}{n}\right), \quad s = 1, \ldots, n, \,\, \xi \in [0,\pi]. 
\end{equation*}
\begin{pro}\label{Prop1}
Let $n \geq 1$ be an integer. Then, for any integer $N\geq0$, we have the following identity
\begin{equation}\label{IdGeg}
\sum_{\substack{m,j \geq 0 \\ N = 2m+nj}}\frac{n(j+k)\Gamma(nk)}{2^{2m+nj}m!\Gamma(n(j+k)+m+1)}C_j^{(k)}(\cos \xi) = \sum_{\substack{j_1, \ldots, j_n \geq 0 \\ j_1+\cdots j_n = N}}  
\frac{(k)_{j_1}\ldots (k)_{j_{n}}}{(nk)_{N}} \frac{b_1^{j_1}}{j_1!}\cdots  \frac{b_{n}^{j_{n}}}{j_{n}!}.
\end{equation}
\end{pro}
A remarkable feature of this identity is that for any $n \geq 2$, its LHS depends on $\xi$ while its RHS does so on the angles $\theta_s$. This strange fact was already present in the variational formula in \cite{Demni1} valid for integer $k \geq 1$, where it comes from the elementary identity
\begin{equation*}
\frac{1}{n}\sum_{s=1}^ne^{2i\pi js/n} = 1
\end{equation*}
if $j \equiv 0[n]$ and vanishes otherwise. It is also present in the Laplace transform \eqref{LapTra}, and can be justified by the formulas proved in the following lemma.

\begin{lem}\label{Lem1}
Let $n\geq 2$ be an integer and let, 
\begin{equation*}
e_m(b_1, \dots, b_n) := \sum_{1 \leq j_1 < \ldots < j_m \leq n} b_{j_1} \ldots b_{j_m}, \quad 1 \leq m \leq n, 
\end{equation*}
be the $m$-th elementary symmetric functions in the variable $(b_1, \ldots, b_n)$. For $m=0$, we set $e_0 :=1$. 
\begin{itemize}
\item If $n$ is odd, then
\begin{equation}\label{El1}
e_n(b_1, \ldots, b_n) = \frac{1}{2^{n-1}} \cos \xi.
\end{equation} 
\item If $n$ is even, then 
\begin{equation}\label{El2}
e_n(b_1, \ldots, b_n) =  \frac{1}{2^{n-1}} \left((-1)^{n/2} - \cos \xi \right). 
\end{equation}
\item For any $n \geq 1$, $e_m(b_1, \dots, b_n) = 0$ for all odd integers $1 \leq m \leq n$, while, 
\begin{equation*}
e_m(b_1,\ldots, b_n) = n (-1)^j \frac{(n-j-1)!}{2^{2j}j!(n-2j)!}, \quad m = 2j, \, \, 0 \leq j < [n/2].
\end{equation*}
\end{itemize}
\end{lem}

Formulas \eqref{El1} and \eqref{El2} are from \cite[table 1.39]{Gra-Ryz}. The rest we postpone to after the proof of Proposition \ref{Prop1} which is the central result of this paper. 
\begin{proof}[Proof of the proposition]
Denote $D_{n,k}(\xi)$ the LHS of \eqref{IdGeg} and multiplying it by $(nk)_Nz^N$ for $|z| < 1, N \geq 0$. Summing the resulting expression over $N \geq 0$, we get 
\begin{equation*}
\sum_{N \geq 0} \Gamma(nk+N) D_{n,k}(\xi) z^N = \sum_{m,j \geq 0}\frac{n(j+k)\Gamma(n(j+k)+2m)}{m!\Gamma(n(j+k)+m+1)}C_j^{k}(\cos \xi) \frac{z^{2m+nj}}{2^{2m+nj}}.
\end{equation*}
Using the following dimidiation formula for the Pochhammer symbol
\[
(x)_{2l}=2^{2l}\Bigl(\frac{x}{2}\Bigr)_{l}\Bigl(\frac{1+x}{2}\Bigr)_{l},
\]
we thus obtain
\begin{equation*}
\sum_{N \geq 0} \Gamma(nk+N) D_{n,k}(\xi) z^N    =  \sum_{j \geq 0} \frac{z^{nj}}{2^{nj}}C_j^{(k)}(\cos \xi) {}_2F_1\left(\frac{nj+nk}{2}, \frac{nj+nk+1}{2}, nj +nk+1; z^2\right). 
\end{equation*}
Thanks to the following formula (see for instance \cite{Erd1}, p.101)
\begin{equation*}
{}_2F_1\left(\frac{nj+nk}{2}, \frac{nj+nk+1}{2}, nj +nk+1; z\right) = \frac{2^{nk+nj}}{(1+\sqrt{1-z})^{nk+nj}}, 
\end{equation*}
valid for $z \in \mathbb{C} \setminus [1,\infty[$, it follows that
\begin{equation*}
\sum_{N \geq 0} \Gamma(nk+N) D_{n,k}(\xi) z^N   = \frac{2^{nk}}{(1+\sqrt{1-z^2})^{nk}}  \sum_{j \geq 0} \frac{z^{nj}}{(1+\sqrt{1-z^2})^{nj}}C_j^{(k)}(\cos \xi), 
\end{equation*}
and, with the use of the generating function for Gegenbauer polynomials
\begin{equation*}
\sum_{j \geq 0} C_j^{(k)}(\cos \xi) w^j = \frac{1}{(1-2w \cos \xi + w^2)^k}, \quad |w| < 1,
\end{equation*}
we are led to
\begin{equation}\label{G1}
\sum_{N \geq 0} \Gamma(nk+N) D_{n,k}(\xi) z^N  = \frac{2^{nk}}{\Bigl((1+\sqrt{1-z^2})^n -2z^n \cos \xi + z^{2n}/(1+\sqrt{1-z^2})^n\Bigr)^k} .
\end{equation}
Moreover, 
\begin{align*}
(1+\sqrt{1-z^2})^n + \frac{z^{2n}}{(1+\sqrt{1-z^2})^n} &= (1+\sqrt{1-z^2})^n + (1-\sqrt{1-z^2})^n
 \\& = 2 \sum_{j=0}^{[n/2]} \binom{n}{2j}(1-z^2)^j = 2z^nT_n\left(\frac{1}{z}\right) ,
 \end{align*}
where $T_n$ is the $n$-th Tchebycheff polynomials (\cite{Erd2}). Now, consider the polynomial equation 
\begin{equation}\label{PolEq}
2z^nT_n\left(\frac{1}{z}\right) - 2\cos(\xi) z^n = 0. 
\end{equation}
If $n$ is odd, then $z^nT_n(z)$ has degree $n-1$ so that \eqref{PolEq} has $n$ roots given 
\begin{equation*}
\frac{1}{b_s}, \quad 1 \leq s \leq n.
\end{equation*}
As a result, \eqref{El1} entails 
\begin{align}\label{G2}
2z^nT_n\left(\frac{1}{z}\right) - 2\cos(\xi) z^n = 2\cos \xi\prod_{s=1}^n \Bigl(\frac{1}{b_s} - z\Bigr)  = 2^n\prod_{s=1}^n \Bigl(1 - b_s z\Bigr).
\end{align}
Otherwise, if $n$ is even then the leading term of $z^nT_n(1/z)$ is $(-1)^{n/2}z^n$ and together with \eqref{El2} lead similarly to
\begin{equation}\label{G3}
2z^nT_n\left(\frac{1}{z}\right) - 2\cos(\xi) z^n  = 2 \left[(-1)^{n/2} - \cos(\xi)\right]\prod_{s=1}^n \Bigl(\frac{1}{b_s} - z\Bigr) = 2^n\prod_{s=1}^n \Bigl(1 - b_s z\Bigr).
\end{equation}
Finally, multiplying the RHS of \eqref{IdGeg} by $z^N, N \geq 0$ and summing over $N$, the generalized binomial Theorem entails
\begin{align}\label{G4}
\sum_{N \geq 0} \left\{\sum_{\substack{j_1, \ldots, j_n \geq 0 \\ j_1+\cdots j_n = N}}  (k)_{j_1}\ldots (k)_{j_{n}} \frac{b_1^{j_1}}{j_1!}\cdots  \frac{b_{n}^{j_{n}}}{j_{n}!}\right\}z^N & = \prod_{s=1}^n \left\{\sum_{j_s \geq 0} \frac{(k)_{j_s}}{j_s!} b_s^{j_s}z^{j_s}\right\}
=  \prod_{s=1}^n\frac{1}{(1 - b_s z)^k}.
\end{align}
Gathering \eqref{G1}, \eqref{G2}, \eqref{G3} and \eqref{G4}, we are done.  
\end{proof}

Now, we shall prove Lemma \ref{Lem1}. 
\begin{proof}[Proof of the lemma]
It remains to prove the last statement of the lemma. To this end, we appeal to Newton identities (\cite[page 23]{MD}):
\begin{equation}\label{Newton}
m e_m(b_1, \ldots, b_n) = \sum_{j=1}^m (-1)^{j-1} e_{m-j}(b_1, \ldots, b_n)p_j(b_1, \ldots, b_n) 
\end{equation}
where 
\begin{equation*}
p_j(b_1, \ldots, b_n) := \sum_{s=1}^n b_s^j, \quad j \geq 0,
\end{equation*}
is the $j$-th Newton sum. Indeed, linearization formulas are available for the Newton sums, namely (\cite[table 1.32]{Gra-Ryz}):
\begin{eqnarray*}
\cos^{2j}(\theta) & = &\frac{1}{2^{2j}} \Biggl(\binom{2j}{j} + 2 \sum_{s=1}^j \binom{2j}{j-s} \cos(2s\theta) \Biggr), \\ 
\cos^{2j+1}(\theta) &=& \frac{1}{2^{2j}} \sum_{s=0}^j\binom{2j+1}{j-s}\cos\bigl((2s+1)\theta\bigr).
\end{eqnarray*}
Varying $\theta \in \bigl\{(\xi+2\pi s)/n, 1 \leq s \leq n\bigr\}$ in the first linearization formula, we see that for any $0 \leq j \leq [m/2]$, 
\begin{equation*}
p_{2j}(b_1, \ldots, b_n) = \frac{1}{2^{2j}}\Biggl(n\binom{2j}{j} + 2 \sum_{s=1}^n\sum_{l=1}^j \binom{2j}{j-l} \cos\left(2l\frac{(\xi + 2s \pi)}{n}\right) \Biggr).
\end{equation*}
Moreover, it is straightforward that for any integer $1 \leq m \leq n-1$ and any $1 \leq l \leq j \leq [m/2]$, 
\begin{equation*}
\sum_{s=1}^n  \cos\left(2l\frac{(\xi + 2s \pi)}{n}\right) = \Re\left(e^{2il\xi/n} \sum_{s=1}^ne^{4ils\pi /n} \right) = 0
\end{equation*}
so that 
\begin{equation*}
p_{2j}(b_1, \ldots, b_n) = \frac{n}{2^{2j}}\binom{2j}{j}.
\end{equation*}
The same reasoning shows that $p_{2j+1}(b_1, \ldots, b_n) = 0$ and an induction on $m$ shows that all elementary symmetric functions $e_m(b_1, \ldots, b_n), 1 \leq m \leq n-1$, do not depend on $\xi$. As a result, 
\begin{align}
e_m(b_1, \ldots, b_n) & = \sum_{1 \leq j_1 < \cdots < j_m \leq n} \cos\left(\frac{2j_1 \pi}{n}\right)  \ldots  \cos\left(\frac{2j_m \pi}{n}\right) \nonumber
\\& = \frac{1}{m}\sum_{j=1}^{[m/2]} (-1)^{2j-1} e_{m-2j}(b_1, \ldots, b_n)p_{2j}(b_1, \ldots, b_n) \nonumber
\\& = \frac{n}{m} \sum_{j=1}^{[m/2]} (-1)^{2j-1} e_{m-2j}(b_1, \ldots, b_n)\frac{1}{2^{2j}}\binom{2j}{j}  \label{Newton1}. 
\end{align}
From \eqref{Newton1}, we readily see by induction that $e_m(b_1,\ldots, b_n) = 0$ if $m$ is odd since $e_1(b_1,\ldots, b_n) = 0$. Otherwise, recall that the $n$-th reverse Tchebycheff polynomial admits the following expansion (\cite[page 185]{Erd2}): 
\begin{equation*}
z^nT_n\left(\frac{1}{z}\right)  = \frac{n}{2} \sum_{j=0}^{[n/2]}(-1)^j \frac{2^{n-2j}(n-j-1)!}{j!(n-2j)!} z^{2j}.
\end{equation*}
Then, for any even integer $m = 2j, 0 \leq j \leq [n/2]$, the expression of $e_m(b_1, \dots, b_n)$ follows from \eqref{G2} and \eqref{G3}. 
\end{proof}

\subsection*{Representation of Gegenbauer polynomials}
When $N=1$, both sides of \eqref{IdGeg} vanish unless $n=1$. Indeed, it is obvious that there is no couple $(m,j) \in \mathbb{N}$ such that $1=2m+nj$ when $n \geq 2$, while the equation 
\begin{equation*}
j_1+\cdots j_n = 1, \quad j_1, \ldots, j_n \geq 0,
\end{equation*}
has exactly $n$ solutions: all but one of these integers vanish in which case the right hand side of \eqref{IdGeg} reduces to 
\begin{equation*}
\frac{1}{n} \sum_{s=1}^n b_s = 0. 
\end{equation*}
If $n=1$, then the only decomposition of $1=2m+j$ corresponds to $(m,j) = (0,1)$ so that \eqref{IdGeg} obviously holds. Now, if $n=2$, then the left hand side vanishes unless $N$ is even. Moreover, if $N = 2M, M \geq 0,$ is even then \eqref{IdGeg} reads
\begin{align*}
\sum_{j=0}^M \frac{2(j+k)\Gamma(2k+2M)}{(M-j)!\Gamma(j+2k+M+1)}C_j^{(k)}(\cos \xi) & = 2^{2M} \frac{(k)_M}{M!} \cos^{2M}(\xi/2)
 \end{align*}
or equivalently (after using the duplication formula for the Gamma function)
\begin{equation*}
\frac{1}{\Gamma(2k)} \left(\frac{1+x}{2}\right)^M = (k+(1/2))_M M!\sum_{j=0}^M \frac{2(j+k)}{(M-j)!\Gamma(j+2k+M+1)}C_j^{(k)}(x), \quad x \in [-1,1].
\end{equation*}
This expansion is an instance of a more general one (see for instance \cite[page 213]{Erd2} or \cite{Koe-Koe}) and is also the inverse (in the sense of composition) of \eqref{ExpGeg}. In this direction, we can extend it to even integers $n = 2q, q \geq 1,$ as follows. If $N = nM, M \geq 0,$ then \eqref{IdGeg} takes the form:
\begin{equation*}
\sum_{j=0}^M \frac{2q(j+k)\Gamma(2qk+2qM)}{(q(M-j))!\Gamma(qj+ qM + 2qk+1)}C_j^{(k)}(\cos \xi) = 2^{2qM}
\sum_{\substack{j_1, \ldots, j_{2q} \geq 0 \\ j_1+\cdots+ j_{2q} = 2qM}}  (k)_{j_1}\ldots (k)_{j_{2q}} \frac{b_1^{j_1}}{j_1!}\cdots  \frac{b_{2q}^{j_{2q}}}{j_{2q}!}.
\end{equation*}
This is an infinite lower-triangular system in the variable $(C_j^{(k)})_{j \geq 0}$ whose diagonal entries $j=M$ equal one. Therefore, it is invertible and one gets after inversion the following corollary of Proposition \ref{Prop1}. 

\begin{cor}
Let $n = 2q, q \geq 1,$ be an even integer and let $k > 0, \xi \in [0,\pi]$ be reals. Then, for any integer $M \geq 0$, there exist real coefficients $(a_j^{(q,k,M)})_{j=0}^M$ such that $a_M^{(q,k,M)} = 1$ and 
\begin{equation*}
C_M^{(k)}(\cos \xi) = \sum_{j=0}^M a_j^{(q,k,M)}
\sum_{\substack{j_1, \ldots, j_{2q} \geq 0 \\ j_1+\cdots+ j_{2q} = 2qj}}  (k)_{j_1}\ldots (k)_{j_{2q}} \frac{b_1^{j_1}}{j_1!}\cdots  \frac{b_{2q}^{j_{2q}}}{j_{2q}!}.
\end{equation*}
\end{cor}
\begin{rems}
\begin{enumerate}
\item It would be interesting to find a compact formula for $(a_j^{(q,k,M)})_{j \geq 0}$ when $q \geq 2$.  
\item When $n$ is odd and $N = nM$, the summation in the left hand side of \eqref{IdGeg} contains either even or odd Gegenbauer polynomials according to whether $M$ is even or odd respectively. Consequently, one gets in both cases analogous representations of $(C_M^{(k)})_{M \geq 0}$. \\
\item In \cite{Dil}, the following representation of Gegenbauer polynomials was derived for integer  $k \geq 1$: 
\begin{equation*}
C_j^{(k)}(\cos \xi) = 2^j \sum_{1 \leq m_1 < \ldots < m_j \leq j+k-1} \,\, \prod_{s=1}^j \left(\cos (\xi) + \cos\left(\frac{m_s \pi}{j+k}\right)\right).
\end{equation*}
It was deduced from a representation of generalized Fibonacci polynomials through elementary symmetric functions (up to a variable change).  
\end{enumerate}
\end{rems}

\section{Concluding remarks}
De Bie et al. formula expressing $F_{n,k}(R, \xi)$ through Horn's function $\Phi_2^{(n)}$ led to the interesting identity \eqref{IdGeg}. The latter has been recently used by the secondly named author in order to solve important problems related to the first hitting time of the boundary of a dihedral wedge by a radial Dunkl process (\cite{Demni2}). On the other hand, the obtained Horn's series is a major step toward the conjectured Laplace-type representation of the generalized Bessel function associated with dihedral groups. 


\begin{thebibliography}{99}
\bibitem{AAR}\emph{G. E. Andrews, R. Askey, R. Roy}. Special functions. {\it Cambridge University Press}. 1999.
\bibitem{BKO}\emph{S. Ben Said, T. Kobayashi, B. Orsted}. Laguerre semigroup and Dunkl operators. {\it Compos. Math}. {\bf 148}, (2012), no. 4, 1265-1336.
\bibitem{CDBL}\emph{D. Constales, H. De Bie, P. Lian}. Explicit formulas for the Dunkl dihedral kenel and the $(\kappa, a)$-generalized Fourier kernel. {\it https://arxiv.org/abs/1610.00098}. 
\bibitem{DeB}\emph{H. De Bie}. The kernel of the radially deformed Fourier transform. {\it Integral Transforms Spec. Funct}. {\bf 24} (2013), no. 12, 1000-1008.
\bibitem{Demni1}\emph{N. Demni}. Radon Transform on spheres and generalized Bessel function associated with dihedral groups. {\it J. Lie Theory},  {\bf 22} (2012), no. 1, 81-91.
\bibitem{Demni2}\emph{N. Demni}. Reciprocal of the First hitting time of the boundary of dihedral wedges by a radial Dunkl process. {\it Submitted}. 
\bibitem{Dil}\emph{K. Dilcher}. A generalization of Fibonacci polynomials and a representation of Gegenbauer polynomials of integer order. {\it Fibonacci Quart}. {\bf 25} (1987), no. 4, 300-303. 
\bibitem{Dij-Koo}\emph{A. Dijksma, T. H. Koornwinder}. Spherical Harmonics and the product of two Jacobi polynomials. {\it Indag. Math.} {\bf 33}, 1971, 191-196.
\bibitem{Erd}\emph{A. Erdelyi, W. Magnus, F. Oberhettinger, F. G. Tricomi}. Tables of Integral Transforms. {\bf Vol. I.} McGraw-Hill Book Company, Inc., New York-Toronto-London, 1954.
\bibitem{Erd1}\emph{A. Erdelyi, W. Magnus, F. Oberhettinger, F. G. Tricomi}. Higher Transcendental Functions. {\bf Vol. I.} McGraw-Hill Book Company, Inc., New York-Toronto-London, 1954.
\bibitem{Erd2}\emph{A. Erdelyi, W. Magnus, F. Oberhettinger, F. G. Tricomi}. Higher Transcendental Functions. {\bf Vol. II.} McGraw-Hill Book Company, Inc., New York-Toronto-London, 1954.
\bibitem{Gra-Ryz} \emph{I.S. Gradshteyn, I.M. Ryzhik}, Table of integrals, series and products, {\it 5th ed., Academic Press, Boston, MA}, 1994.
\bibitem{Kob-Man}\emph{T. Kobayashi, G. Mano}. The inversion formula and holomorphic extension of the minimal representation of the conformal group. {\it Harmonic analysis, group representations, automorphic forms and invariant theory}, 151-208, 
Lect. Notes Ser. Inst. Math. Sci. Natl. Univ. Singap., 12, World Sci. Publ., Hackensack, NJ, 2007. 
\bibitem{Koe-Koe}\emph{J. Koekoek, R. Koekoek}. The Jacobi inversion formula. {\it Complex Variables}, {\bf 39}, (1999), 1-18.
\bibitem{MD}\emph{I. G. MacDonald}. Symmetric Functions and Hall Polynomials. {\it Second edition, Mathematical Monographs, Oxford}. 1995.
\bibitem{Wat}\emph{G. N. Watson}. A treatise on the theory of Bessel functions. {\it Second Edition. Cambridge Mathematical Library edition}, 1962.  
\bibitem{Xu}\emph{Y. Xu}. A product formula for Jacobi polynomials. {\it Special functions (Hong Kong, 1999), 423-430, World Sci. Publ., River Edge, NJ}, 2000. 
\end{thebibliography}
\end{document}